\documentclass{sampta}
\usepackage{graphicx}
\usepackage{type1cm}
\usepackage{eso-pic}
\usepackage{color}
\usepackage{amsmath,amssymb,amsthm}
\usepackage{enumerate}

\newtheorem{thm}{Theorem}
\newtheorem{prop}[thm]{Proposition}
\newtheorem*{RP}{Piecewise $D$-finite Reconstruction Problem}

\DeclareMathOperator{\diffop}{\mathcal{E}} 
\DeclareMathOperator{\shift}{E} \DeclareMathOperator{\id}{I}
\DeclareMathOperator{\Op}{\mathfrak{D}} 
\DeclareMathOperator{\newop}{\ensuremath{\widehat{\Op}}}

\newcommand{\nullsp}[1]{\ensuremath{\mathcal{N}_{#1}}} 
\renewcommand{\vec}[1]{\ensuremath{\mathbf{#1}}} 
\newcommand{\f}[1]{\mbox{$#1$}}
\newcommand{\isdef}{\ensuremath{\stackrel{\text{def}}{=}}}
\newcommand{\reals}{\ensuremath{\mathbb{R}}}
\newcommand{\naturals}{\ensuremath{\mathbb{N}}}
\newcommand{\step}{\ensuremath{\mathcal{H}}}
\newcommand{\fn}{\ensuremath{g}} 
\newcommand{\n}[1]{\ensuremath{\widetilde{\fn_{#1}}}} 
\newcommand{\theorder}{\ensuremath{N}} 
\newcommand{\polymaxorder}{\ensuremath{k}} 
\newcommand{\np}{\ensuremath{\mathcal{K}}} 
\newcommand{\coeff}{\ensuremath{p}} 
\newcommand{\polycoeff}{\ensuremath{a}} 

\newcommand{\nm}{\ensuremath{M}} 
\newcommand{\nmH}{\ensuremath{\widehat \nm}}
\newcommand{\nmC}{\ensuremath{\widetilde \nm}}
\newcommand{\stepdefinition}{\ensuremath{\begin{cases}0 & x<0\\ 1 & x \geq 0\end{cases}}} 

\newcommand{\nc}{\newcommand}
\nc{\be}{\begin{equation}}
\nc{\ee}{\end{equation}}
\nc{\Rr}{{\mathbb R}}
\nc{\Nn}{{\mathbb N}}
\nc{\Ff}{{\mathcal F}}
\def\phi{\varphi}
\title{{\titleitfnt An ``algebraic'' reconstruction of piecewise-smooth functions
from integral measurements}}

\author{
Dima Batenkov, Niv Sarig, Yosef Yomdin\\[12pt]
\addrfnt Department of Mathematics, Weizmann institute of science, Rehovot, Israel.\\
\email{\{dima.batenkov, niv.sarig, yosef.yomdin\}@weizmann.ac.il}
}

\begin{document}

\maketitle

\section{Introduction}

This paper presents some results on a well-known problem in
Algebraic Signal Sampling and in other areas of applied
mathematics: reconstruction of piecewise-smooth functions from
their integral measurements (like moments, Fourier coefficients,
Radon transform, etc.). Our results concern reconstruction (from
the moments) of signals in two specific classes: linear
combinations of shifts of a given function, and ``piecewise
$D$-finite functions" which satisfy on each continuity interval a
linear differential equation with polynomial coefficients.

Let us start with some general remarks and a conjecture. It is
well known that the error in the best approximation of a
$C^k$-function $f$ by an $N$-th degree Fourier polynomial is of
order ${C\over {N^k}}$. The same holds for algebraic polynomial
approximation and for other basic approximation tools. However,
for $f$ with singularities, in particular, with discontinuities,
the error is much larger: its order is only ${C\over {\sqrt N}}$.
Considering the so-called Kolmogorow $N$-width of families of
signals with moving discontinuities one can show that {\it any
linear approximation method provides the same order of error, if
we do not fix a priori the discontinuities' position} (see
\cite{Ett.Sar.Yom}, Theorem 2.10). Another manifestation of the
same problem is the {\it ``Gibbs effect" - a relatively strong
oscillation of the approximating function near the
discontinuities}. Practically important signals usually do have
discontinuities, so the above feature of linear representation
methods presents a serious problem in signal reconstruction. In
particular, it visibly appears near the edges of images compressed
by JPEG, as well as in the noise and low resolution of the CT and
MRI images.

Recent non-linear reconstruction methods, in particular,
Compressed Sensing (\cite{Can,Don}) and Algebraic Sampling
(\cite{Vet1,Vet2,Vet3,Mil,Put}), address this problem in many
cases. Both approaches appeal to an a priori information on the
character of the signals to be reconstructed, assuming their
``simplicity" in one or another sense. Compressed sensing assumes
only a sparse representation in a certain (wavelets) basis, and
thus it presents a rather general and ``universal" approach.
Algebraic Sampling usually requires more specific a priori
assumptions on the structure of the signals, but it promises a
better reconstruction accuracy. In fact, we believe that
ultimately the Algebraic Sampling approach has a potential to
reconstruct ``simple signals with singularities" as good as smooth
ones. In particular, the results of \cite{Eck,Kve,Tad1,Tad2,Vet3}
strongly support (also apparently do not accurately formulate and
prove) the following conjecture:

{\it There is a non-linear algebraic procedure reconstructing any
signal in a class of piecewise $C^k$-functions (of one or several
variables) from its first $N$ Fourier coefficients, with the
overall accuracy of order ${C\over {N^k}}$. This includes the
discontinuities' positions, as well as the smooth pieces over the
continuity domains.}

At present there are many approaches available to a robust
detection of discontinuities from Fourier data (see
\cite{Tad1,Eck,Kve} and references therein). The remaining problem
seems to be an accurate estimate of the accuracy of the solution
of the nonlinear systems arising. Our results below can be
considered, in particular, as a step in this direction. On the
other hand, they have been motivated by the results in
\cite{Vet1,Vet2,Vet3}, and in \cite{Put, Mil}.

\section{Linear combinations of shifts of a given function}

Reconstruction of this class of signals from sampling has been
described in \cite{Vet1,Vet2}. We study a rather similar problem
of reconstruction from the moments. Our method is based on the
following approach: we construct convolution kernels dual to the
monomials. Applying these kernels, we get a Prony-type system of
equations on the shifts and amplitudes.

Let us restate a general reconstruction problem, as it appears in
our specific setting. We want to reconstruct signals of the form
\be\label{eq:the_model} F(x)=\sum_{i=1}^N\sum_{j,l}a_{i,j,l}
f_i^{(l)}(x+x^j) \ee where the $f_i$'s  are known functions of
$x=(x_1,\ldots,x_d)$, and the form (\ref{eq:the_model}) of the
signal is known a priori. The parameters $a_{i,j,l}, \
x^j=(x^j_1,\dots,x^j_d)$ are to be found from a finite number of
``measurements", i.e. of linear (usually integral) functionals
like polynomial moments, Fourier moments, shifted kernels,
evaluation over some grid and more.

In this paper we consider only linear combinations of shifts of
one known function $f$ (although the method of ``convolution dual"
can be extended to several shifted functions and their derivatives
- see \cite{Sar.Yom2}). First we consider general integral
``measurements" and then restrict ourselves to the moments and
Fourier coefficients. In what follows $x=(x_1,\dots,x_d),t=(t_1,\ldots,t_d)$, $j$ \
is a scalar index, while $k=(k_1,\dots,k_d), \ i=(i_1,\dots,i_d)$
and $n=(n_1,\dots,n_d)$ are multi-indices. Partial ordering of
multi-indices is given by $k \leq k' \Leftrightarrow k_p \leq
k'_p, \ p=1,\dots,d.$ So we have \be\label{eq:model_translations}
F(x)=\sum_{j=1}^s a_j f(x+x^j). \ee Let the measurements
$\mu_k(F)$ be given by $\mu_k(F)=\int F(t) \phi_k(t)dt,$ for a
certain (multi)-sequence of functions $\phi_k(t), \ k\geq
0=(0,\dots,0)$.

Given $f$ and $\phi=\{\phi_k(t)\}, \ k\geq 0$ we now try to find
certain ``triangular" linear combinations
\be\label{eq:triangular_system} \psi_k(t)=\sum_{0\leq i\leq k}
C_{i,k}\phi_i(t) \ee forming, in a sense, some ``$f$-convolution
dual" functions (similar to a bi-orthogonal set of function) with
respect to the system $\phi_k(t)$. More accurately, we require
that \be\label{eq:dual_convolution} \int f(t+x)\psi_k(t) =
\phi_k(x).\ee We shall call a sequence $\psi = \{\psi_k(t)\}$
satisfying (\ref{eq:triangular_system}), \
(\ref{eq:dual_convolution}) \ $f$ - convolution dual to $\phi$.
Below we find convolution dual systems to the usual and
exponential monomials.

We consider a general problem of finding convolution dual
sequences to a given sequence of measurements as an important step
in the reconstruction problem. Notice that it can be generalized
by dropping the requirement of a specific representation
(\ref{eq:triangular_system}): $\psi_k(t)=\sum_{i=0}^k
C_{i,k}\phi_i(t)$. Instead we can require only that $\int
f(t)\psi_k(t)$ be expressible in terms of the measurements
sequence $\mu_k$. Also $\phi_k$ in (\ref{eq:dual_convolution}) can
be replaced by another a priori chosen sequence $\eta_k$. This
problem leads, in particular, to certain functional equations,
satisfied by polynomials and exponents (as well as exponential
polynomials and some kinds of elliptic functions).

Now we have the following result:

\begin{thm}\label{thm:conv_dual} Let a sequence $\psi = \psi_k(t)$ be $f$-convolution
dual to $\phi$. Define $M_k$ by $M_k=\sum_{0\leq i\leq k}
C_{i,k}\mu_i.$ Then the parameters $a_j$ and $x^j$ in
(\ref{eq:model_translations}) satisfy the following system of
equations (``generalized Prony system"):
\be\label{eq:gnrlzd_prony} \sum_{j=1}^s a_j\phi_k(x^j)=M_k, \ \
k=0, \dots.\ee
\end{thm} {\bf Proof} We have $M_k=\sum_{0\leq i\leq k} C_{i,k}\mu_i = \int
F(t)\sum_{0\leq i\leq k} C_{i,k}\phi_i(t) dt = \int
F(t)\psi_k(t)=\sum_{j=1}^s a_j \int f(t+x^j)\psi_k(t) dt =
\sum_{j=1}^s a_j\phi_k(x^j).$

In specific examples we can find the minimal number of equations
in (\ref{eq:gnrlzd_prony}) necessary to uniquely reconstruct the parameters $a_j$ and
$x^j$ in (\ref{eq:model_translations}).

\subsection{Reconstruction from moments}\label{section_poly} We
are given a finite number of moments of a signal $F$ as in
(\ref{eq:model_translations}) in the form \be m_n=\int F(t)t^ndt.
\ee So here $\phi_n(x)=x^{n_1}_1 \cdots x^{n_d}_d$ for each
multi-index $n=(n_1,\dots,n_d)$. We look for the dual functions
$\psi_n$ satisfying the convolution equation
\be\label{eq:conv_poly} \int f(t+x)\psi_n(t)dt=x^n \ee for each
multi-index $n$. To solve this equation we apply Fourier transform
to both sides of (\ref{eq:conv_poly}). Assuming that $\hat
f(\omega)\in C^\infty(\Rr^d),\hat f(0)\neq0$ we find (see
\cite{Sar.Yom2}) that there is a unique solution to
(\ref{eq:conv_poly}) provided by \nc{\pa}{\partial}
\be\label{eq:dual_poly} \phi_n(x)=\sum_{k\leq n}C_{n,k}x^k, \ee
where
\[
C_{n,k}=\frac1{(\sqrt{2\pi})^d}{n\choose k} (-i)^{n+k} \left[
\left.\frac{\pa^{n-k}}{\pa\omega^{n-k}}
\right|_{\omega=0}\frac1{\hat f(\omega)}\right].
\]
So we set the generalized polynomial moments as
\be\label{eq:gnrlzd_moments_poly} M_n=\sum_{k\leq n}C_{n,k}m_k \ee
and obtain, as in Theorem \ref{thm:conv_dual}, the following
system of equations: \be\label{eaution_prony_poly}
\sum^s_{j=1}a_j(x^j)^n = M_n, \ n \geq 0.\ee This system can be
solved explicitly in a standard way (see, for example,
\cite{Nik.Sor,Vet1,Sar.Yom1}). In one-dimensional case it goes as
follows (see \cite{Nik.Sor}): from (\ref{eaution_prony_poly}) we
get that for $z=(z_1,\dots,z_d)$ the generalized moments
generating function \be\label{eq:genrating_function} I(z)
=\sum_{n\in \Nn^d}M_nz^n =\sum_{j=1}^s
a_j\prod_{l=1}^d\frac{1}{1-x^j_lz_l} \ee is a rational function.
Hence its Taylor coefficients satisfy linear recurrence relation,
which can be reconstructed through a linear system with the
Hankel-type matrix formed by an appropriate number of the moments
$M_n$'s. This is, essentially, a procedure of the diagonal Pad\'e
approximation for $I(z)$ (see \cite{Nik.Sor}). The parameters
$a_j, x^j$ are finally reconstructed as the poles and the residues
of $I(z)$. For several variables the solution procedure is
similar.
\\In one dimensional case with the derivatives $f^{(l)}$ included  we have
\be\label{eq:model_trans_der}
F(x)=\sum^s_{j=1}\sum^r_{l=0}a_{j,l}f^{(l)}(x+x^j).\ee The
corresponding moment-generating function in this case takes the
form \be\label{eq:generating_der} I(z)
=\sum_{j=1}^s\sum_{l=0}^r\sum_{q=0}^{l}{l\choose q} \frac
{(-1)^{q+l}a_{j,l}/(x^j)^{l}}{(1-x^jz)^{q+1}}. \ee which is still
a rational function (d-dimensional case with derivatives is
similar). We would like to stress that in this case the dual
polynomials $\psi_k$ are not changed and they are given as in
(\ref{eq:dual_poly}). Therefore also the formula for the
generalized moments $M_n$ is the same as in
(\ref{eq:gnrlzd_moments_poly}).

\subsection{Fourier case}\label{section_fourier}
In the same manner as in section \ref{section_poly} we now choose
$\phi_k(x)=e^{ikx}$. We get immediately $\psi_k(x)=\frac1{\hat
f(k)}e^{-ikx}$. Indeed,
\[
\int f(t+x)\psi_k(t)dt=\int f(t+x)\frac1{\hat f(k)}e^{ikt}dt=
\]
\be\label{eq:fourier} \frac{\hat f(k)}{\hat
f(k)}e^{-ikx}=\phi_{-k}(x). \ee Here the triangular system of
equations (\ref{eq:triangular_system}) is actually not triangular any more but still
since $\psi_k(x)=\frac1{ \hat f(k)}\phi_{-k}(x)$ we can express the generalized moments through the original ones via $M_k=
\frac1{\hat f(k)}\mu_{-k}[F]$. Now exactly as
before we can find a generalized Prony system in the form
\be\label{eq:prony_fourier} \frac1{\hat
f(k)}\mu_{-k}[F]=M_k=\sum_ja_je^{-ikx_j}=\sum_ja_j\rho_j^{k} \ee where
$\rho_j=e^{-ix_j}$. In this case we get a rational exponential
generating function and we can find its poles and residues on the
unit complex circle as we did in the polynomial case.
\subsection{Further extensions}\label{subsection_extension}
The approach above can be extended in the following directions: 1.
Reconstruction of signals built from several functions or with the addition of dilations also can be
investigated (a perturbation approach where the dilations
are approximately 1 is studied in \cite{Sar.Yom1}). 2. Further study of
``convolution duality" can significantly extend the class of
signals and measurements allowing for a closed - form signal
reconstruction.

\section{Reconstruction of piecewise $D$-finite functions from moments}

Let \f{\fn:[a,b] \to \reals} consist of $\np+1$ ``pieces''
\f{\fn_0,\dotsc \fn_\np} with \f{\np \geq 0} jump points
\[
a=\xi_0 < \xi_1 \dotsc <\xi_{\np} < \xi_{\np+1}=b
\]
Furthermore, let $\fn$ satisfy on each continuity interval some
linear homogeneous differential equation with polynomial
coefficients: \f{\Op \fn_n \equiv 0, \; n=0,\dotsc,\np} where
\begin{equation}\label{eq:operator}
\mathfrak{D} = \sum_{j=0}^\theorder
\biggl(\sum_{i=0}^{\polymaxorder_j} \polycoeff_{i,j} x^i \biggr)
\frac{d^j}{dx^j} \quad (\polycoeff_{ij} \in \reals)
\end{equation}
Each \f{\fn_n} may be therefore written as a linear combination of
functions \f{\{u_i\}_{i=1}^\theorder} which are a basis for the
space \f{\nullsp{\Op} = \{f: \Op f \equiv 0\}}:
\begin{equation}\label{eq:linear-comb-nullsp}
\fn_n(x)=\sum_{i=1}^\theorder \alpha_{i,n} u_i(x), \quad n=0,1,\dotsc,\np
\end{equation}
We term such functions $\fn$ ``piecewise $D$-finite''. Many
real-world signals may be represented as piecewise $D$-finite
functions, in particular: polynomials, trigonometric functions,
algebraic functions.

The sequence \f{\{m_k=m_k(\fn)\}} is given by the usual moments
\[
m_k(\fn) = \int_a^b x^k \fn(x) dx
\]
We subsequently formulate the following
\begin{RP} Given \f{\theorder, \{\polymaxorder_i\}, \np,a,b} and
the moment sequence \f{\{m_k\}} of a piecewise $D$-finite function
$\fn$, reconstruct all the parameters
\f{\{\polycoeff_{i,j}\},\{\xi_i\}, \{\alpha_{i,n}\}}.
\end{RP}
Below we state some results (see \cite{Bat} for detailed proofs)
which provide \emph{explicit algebraic connections} between  the
above parameters and the measurements \f{\{m_k\}}.

The first two theorems assume a single continuity interval
(compare \cite{Kis}).

\newcommand{\vv}[2]{v^{(#1)}_{#2}}
\begin{thm}\label{thm:recurrence1}
Let \f{\np=0} and \f{\Op \fn \equiv 0} with $\Op$ given by
\eqref{eq:operator}. Then the moment sequence \f{\{m_k(\fn)\}}
satisfies a linear recurrence relation
\begin{align}\label{eq:recurrence}
\biggl((\shift-a\id)^{\theorder} (\shift-b\id)^\theorder \cdot
\sum_{j=0}^{\theorder} \sum_{i=0}^{\polymaxorder_j} a_{i,j}
\Pi^{(i,j)}(k,\shift) \biggr) m_k = 0
\end{align}
where $\shift$ is the discrete forward shift operator and
\f{\Pi^{(i,j)}(k,\shift)} are monomials in $\shift$ whose
coefficients are polynomials in $k$:
\f{\Pi^{(i,j)}(k,\shift)=(-1)^j \frac{(i+k)!}{(i+k-j)!}
\shift^{i-j}}.
\end{thm}

\begin{thm}\label{thm:system}
Denote
\begin{align*}
\diffop(\shift) &\isdef (\shift-a\id)^{\theorder}
(\shift-b\id)^\theorder, &
\vv{i,j}{k} &\isdef \bigl(\diffop(\shift) \cdot \Pi^{(i,j)}(k,\shift) \bigr)m_k,\\
h_j(z) &\isdef \sum_{k=0}^\infty \vv{0,j}{k} z^k, &
G_j(x) &\isdef \diffop(x) \frac{d^j}{dx^j} \fn(x)
\end{align*}
Assume the conditions of Theorem \ref{thm:recurrence1}. Then
\begin{enumerate}[(1)]
\item The vector of the coefficients \f{\vec{\polycoeff}=(\polycoeff_{i,j})}
satisfies a linear homogeneous system
\begin{align}\label{eq:systemH}
H \vec{\polycoeff} = \begin{pmatrix}
\vv{0,0}{0} & \vv{1,0}{0} & \dots & \vv{\polymaxorder_{\theorder},\theorder}{0} \\
\vv{0,0}{1} & \vv{1,0}{1} & \dots & \vv{\polymaxorder_{\theorder},\theorder}{1} \\
\vdots & \vdots & \vdots & \vdots \\
\vv{0,0}{\nmH} & \vv{1,0}{\nmH} & \dots & \vv{\polymaxorder_{\theorder},\theorder}{\nmH} \\
\end{pmatrix} \begin{pmatrix}
\polycoeff_{0,0} \\
\polycoeff_{1,0} \\
\vdots \\
\polycoeff_{\polymaxorder_{\theorder},\theorder} \\
\end{pmatrix}= 0
\end{align}
for all \f{\nmH \in \naturals}.
\item \f{\vv{i,j}{k}=m_{i+k}\left( G_j(x) \right)}. Consequently,
$h_j(z)$ is the moment generating function of $G_j(x)$.
\item Denote \f{p_j(x)=\sum_{i=0}^{\polymaxorder_j} \polycoeff_{i,j} x^i}.
Then the functions \f{\Phi=\{1,h_0(z), \dotsc h_{\theorder}(z)\}}
are polynomially dependent:
\begin{align*}
\sum_{j=0}^{\theorder} h_j(z) \bigl( z^{\max \polymaxorder_j} \coeff_j(z^{-1}) \bigr) = Q(z)
\end{align*}
where $Q(z)$ is a polynomial with \f{\deg Q < \max \polymaxorder_j}.
The system of polynomials \f{\{z^{\polymaxorder_j} \coeff_j(z^{-1})\}}
is called the Pad\'{e}-Hermite form for $\Phi$.
\end{enumerate}
\end{thm}
To handle the piecewise case, we represent the jump
discontinuities by the step function \f{\step(x) \isdef
\stepdefinition} and write $\fn$ as a distribution
\begin{align}\label{eq:piecewise-def}
\fn(x) &= \n{0}+\sum_{n=1}^\np \n{n}(x) \step(x-\xi_n)
\end{align}
\begin{thm}\label{thm:piecewise}
Let \f{\np > 0} and let $\fn$ be as in \eqref{eq:piecewise-def}
with operator $\Op$ annihilating every piece $\n{n}$. Then the
operator
\begin{align}\label{eq:newop}
\newop \isdef \biggl\{\prod_{n=1}^\np (x-\xi_i)^\theorder \id \biggr\}\cdot \Op
\end{align}
annihilates the entire $\fn$ as a distribution. Consequently,
conclusions of Theorems \ref{thm:recurrence1} and \ref{thm:system}
hold with $\Op$ replaced by $\newop$ as in \eqref{eq:newop}.
\end{thm}

\begin{prop}
Let \f{\np\geq 0} and \f{\{u_i\}_{i=1}^\theorder} be a basis for
the space \f{\nullsp{\Op}}, where $\Op$ annihilates every piece of
$\fn$. Assume \eqref{eq:linear-comb-nullsp} and denote
\f{c_{i,k}^n=\int_{\xi_n}^{\xi_{n+1}} x^k u_i(x)} for
\f{n=0,\dotsc,\np}. A straightforward computation gives \f{\forall
\nmC \in \naturals}:
\begin{align}\label{eq:particular-solution}
\begin{pmatrix}
c_{1,0}^0 & \dotsc & c_{\theorder,0}^0 & \dotsc & c_{\theorder,0}^\np\\
\vdots & \vdots & \vdots & \vdots  & \vdots\\
c_{1,\nmC}^0 & \dotsc & c_{\theorder,\nmC}^0 & \dotsc & c_{\theorder,\nmC}^\np
\end{pmatrix}
\begin{pmatrix}
 \alpha_{1,0}\\ \vdots \\ \alpha_{\theorder,0}\\
\vdots \\ \alpha_{\theorder,\np}
\end{pmatrix}
=
\begin{pmatrix}
 m_0\\ m_1\\ \vdots\\ m_{\nmC}
\end{pmatrix}
\end{align}
\end{prop}

The above results can be combined as follows to provide a solution
of the Reconstruction Problem:
\begin{enumerate}[(a)]
\item Let \f{\theorder, \{\polymaxorder_i\}, \np,a,b} and \f{\{m_k(\fn)\}}
be given. If \f{\np >0}, replace \f{\Op} (still unknown) with
\f{\newop} according to \eqref{eq:newop}.
\item Build the matrix $H$ as in \eqref{eq:systemH}. Solve \f{H\vec{a}=0}
and obtain the operator \f{\Op^*=\Op_{\vec{a}}} which annihilates $\fn$.
\item If \f{\np>0}, factor out all the common roots of the polynomial
coefficients of $\Op^*$ with multiplicity $\theorder$. These are
the locations of the jump points \f{\{\xi_n\}}. The remaining part
is the operator $\Op^{\dagger}$ which annihilates every $\fn_n$.
\item By now $\Op^\dagger$ and $\{\xi_n\}$ are known. So compute
the basis for \f{\nullsp{\Op^{\dagger}}} and solve \eqref{eq:particular-solution}.
\end{enumerate}

The constants $\nmH$ and $\nmC$ determine the minimal required
size of the corresponding linear systems \eqref{eq:systemH} and
\eqref{eq:particular-solution} in order for all the solutions of
these systems to be also solutions of the original problem. It can
be shown that:

\begin{enumerate}
\item There exists no uniform bound on $\nmH$ without any additional
information on the nature of the solutions. Explicit bounds may be
obtained for simple function classes such as piecewise polynomials
of bounded degrees or real algebraic functions.
\item For every specific $\Op$, an explicit bound \f{\nmC=\nmC(\Op)}
may be computed for the system \eqref{eq:particular-solution}.
\end{enumerate}
The above algorithm has been tested on exact reconstruction of
piecewise polynomials, piecewise sinusoids and rational functions.

\end{document}